\documentclass[a4paper,12pt,reqno]{amsart}
\usepackage[utf8]{inputenc}
\usepackage{csquotes}
\usepackage{amssymb}
\usepackage{amsmath}
\usepackage{makecell}
\usepackage{multirow}
\usepackage{hhline}
\usepackage{amsthm}
\usepackage{amsfonts}
\usepackage{mathtools}
\usepackage{todonotes}
\usepackage{comment}
\usepackage{mathrsfs}
\usepackage{enumitem}
\usepackage[all,cmtip]{xy}
\usepackage{tikz-cd}
\usepackage[hyphens]{url}
\usepackage{hyperref}
\usepackage{cleveref}
\usepackage[style=alphabetic]{biblatex}
\usepackage[british]{babel}
\usepackage[margin=6em]{geometry}
\usepackage{eqparbox}

\usepackage{pbox}

\renewcommand{\O}{\mathcal{O}} 

\newcommand{\set}[1]{\left\lbrace #1 \right\rbrace}

\newcommand{\field}[1]{\mathbb{#1}}  
\newcommand{\Z}{\field{Z}} 

\renewcommand{\P}{\field{P}}
\newcommand{\PP}{\field{P}}

\DeclareMathOperator{\ddiv}{div}

\DeclareMathOperator{\Pic}{Pic}

\DeclareMathOperator{\Sym}{Sym}

\DeclareMathOperator{\genus}{g}

\newtheorem{lemma}{Lemma}
\newtheorem{theorem}[lemma]{Theorem}
\newtheorem{proposition}[lemma]{Proposition}
\newtheorem{corollary}[lemma]{Corollary}

\theoremstyle{definition}
\newtheorem{definition}[lemma]{Definition}

\numberwithin{lemma}{section}
\numberwithin{equation}{section} 
\numberwithin{figure}{section}
\title{Contracting low degree points on curves.}
\author{Maarten Derickx}

\begin{document}

\begin{abstract}
The main result of this article is that all but finitely many points of small enough degree on a curve can be written as a pullback of an even smaller degree point. The main theorem has several corollaries that yield improvements on results of Kadets and Vogt, Khawaja and Siksek, and Vojta under a slightly stronger assumption on the degree of the points.
\end{abstract}

\maketitle

\vspace{-1.5em}

\section{Introduction}

The main result of this article is the following Theorem:

\begin{theorem}\label{thm:main}
Let $X$ and $Y$ be nice\footnote{nice here means smooth, projective, and geometrically integral} curves over a number field $K$, and $f: X \to Y$ be a dominant morphism. Let $d$ be an integer such that 
\begin{align}
   \genus(X) > \deg(f) \genus(Y) + (2d-1)(\deg f-1) \label{eq:cs_ineq}.
\end{align}
Then there is a finite set $T_d$ of closed degree $d$ points such that for all closed points $x \in X \setminus T_d$ of degree $d$ there exists a decomposition $X \stackrel{f_1}{\to} Z \stackrel{f_2}{\to} Y$ of $f$ with $\deg f_1 > 1$ such that $x$ comes from pulling back a closed point of degree $d/\deg f_1$ on $Z$.
\end{theorem}
This theorem has several corollaries, that are related to existing results.

\begin{corollary}\label{cor:vojta}
    Under the same assumptions of \Cref{thm:main} the set 
    $$\set{x \in X \text{ closed of degree } \leq d \mid K(f(x)) = K(x) }$$
    is contained in $\cup_{i=1}^d T_i$ and hence finite.
\end{corollary}

\begin{corollary}\label{cor:kadets_vogt}
Let $X$ be a nice curve over a number field $K$ and let $d$ be an integer such that $\genus(X)> (2d-1)^2$. Then there exist a finite set $S$ consisting of morphisms of curves of the form $f_1: X \to Z$ where $Z$ is a nice curve and $\deg f_1 > 1$, such that all but finitely many closed degree $d$ points on $X$ come from pulling back a degree $d/\deg f_1$ point along some $f_1 \in S$.
\end{corollary}
\begin{proof} Since $X$ has infinitely many closed degree $d$ points there is an $f \in K(X)$ of degree $\leq 2d$ \cite[Prop. 2]{frey1994curves}. 
Apply \Cref{thm:main} with this $f$. The corollary follows since up to isomorphism there are only finitely many decompositions $X \stackrel{f_1}{\to} Z \stackrel{f_2}{\to} Y$ .
\end{proof}
\begin{corollary}[Single Source Theorem]\label{cor:kawhaja_siksek}
Let $X$ be a nice curve over a number field $K$. Let $c$ be an integer such that $X$ has infinitely many primitive points of degree $c$ and $\genus(X)> (2c-1)^2$. Let $d$ be an integer such that $\genus(X)> 2d(c-1)+1$.
\begin{enumerate}
    \item There exists a nice curve $Y$ of over $K$ such that $\# Y(K) = \infty$ and an indecomposable covering $f: X \to Y$ of degree $c$ such that all but finitely many primitive points of degree $c$ are of the form $f^{*}(x)$ for some $x \in Y(K)$.\label{cor:kawhaja_siksek2}
    \item If $d\neq c$, then $X$ has finitely many primitive points of degree $d$.  \label{cor:kawhaja_siksek1}
    \item If $X$ has infinitely many closed points of degree $d$ then $c \mid d$.\label{cor:kawhaja_siksek3}
    \item All but finitely many closed points of degree $d$ are of the form $f^*(z)$ where $z$ is a closed point on $Z$ is of degree $d/c$. \label{cor:kawhaja_siksek4}
\end{enumerate}
\end{corollary}

The first two items in \Cref{cor:kawhaja_siksek} can be seen as generalisation of the recent Single Source Theorem by Kawhaja and Siksek, as explained in \Cref{sec:single_source_comp}. However, the last two items are something new and give interesting results regarding the density degree sets studied by Kadets and Vogt.

\begin{definition}\cite[p. 1]{KadetsVogt}
Let $X$ be a nice curve over a number field $K$. Then the \textit{density degree} set $\delta(X/K)$ is the set of all positive integers $d$ such that $X$ has infinitely many closed points of degree $d$. The \textit{minimum density degree} of $X$ is defined as $\min(\delta(X/K)).$
\end{definition}

With this definition, one can reinterpret \Cref{cor:kawhaja_siksek3} of \Cref{cor:kawhaja_siksek}. Indeed, it shows that the existence of infinitely many primitive points of some small degree $c$, forces the integers in $\delta(X/K)$ that are small enough to be divisible by $c$.

\subsection{Acknowledgements}
I would like to thank Borys Kadets and Isabel Vogt for the useful discussion while writing this paper. I would like to thank Maleeha Khawaja for suggesting some corrections to an earlier version of the paper.
The author was supported by the Croatian Science Foundation under project No. IP-2022-10-5008 during this research.

\section{Comparison with existing results.}
\subsection{Comparing \texorpdfstring{\Cref{cor:vojta}}{Corollary 1.2}}

This corollary is similar to the following result.

\begin{theorem}[{\cite[Prop. 2.3]{SongTucker}}]\label{thm:song-tucker}
Let $X$ and $Y$ be nice curves over a number field $K$, and $f: X \to Y$ be a dominant morphism. Assume that $$\genus(X) -1>(d+\genus(Y)-1)\deg f.$$ Then the following set is finite: $$\set{x \in X \text{ closed of degree } \leq d \mid K(f(x)) = K(x) }.$$ 
\end{theorem}

The above theorem with $\genus(Y)=0$ was already proved by Voijta \cite[Cor 0.3]{VojtaGeneralization}. Note that \cref{eq:cs_ineq} used in \Cref{cor:vojta} can be rewritten as $$\genus(X) -1>(d+\genus(Y)-1)\deg f + d(\deg f-2).$$ So the bound of \Cref{cor:vojta} is worse and differs by $d(\deg f-2)$ from \Cref{thm:song-tucker}. Note that \Cref{cor:vojta} would be better if $\deg f = 1$, but this is incompatible with \cref{eq:cs_ineq}. It also means that in the case $\deg f=2$ the bounds of \Cref{thm:song-tucker} and \Cref{cor:vojta} are the same. The proof of \Cref{cor:vojta} can be seen as an alternative proof of \Cref{thm:song-tucker}, where the height machinery has been hidden in the work of Faltings \cite{FaltingsAV} on the rational points of a subvariety of an abelian variety \cite{FaltingsAV}. 

Although \Cref{thm:song-tucker} gives a better or equal bound in all cases, there are advantages to our proof of \Cref{cor:vojta}. That is, \Cref{thm:main} also gives a geometric explanation for why $K(f(x)) \subsetneq K(x)$ for all but finitely many closed points. Indeed, it shows that all but finitely many closed points of small enough degree are contracted to a closed point of degree $d/\deg f_1$ by some $f_1$ which is a partial factorisation of $f$, this does not immediately follow from \Cref{thm:song-tucker}. It is this improvement that allows for the other two corollaries.

\subsection{Comparing \texorpdfstring{\Cref{cor:kadets_vogt}}{Corollary 1.3}}

\Cref{cor:kadets_vogt} is a generalisation of a result by Kadets and Vogt. Here, it must be noted that \Cref{cor:kadets_vogt} needs a stronger genus condition to make this generalisation.

\begin{theorem}[{\cite[Theorem 1.3]{KadetsVogt}}]\label{thm:kadets_vogt}
Suppose that $X$ is a nice curve over a number field $K$. Let $d$ the minimum density degree of $X$, let $m:=\lceil{d/2}\rceil -1$ and let $\epsilon := 3d-1-6m < 6.$ If 
\begin{align}\genus(X) > \max \left(\frac{d(d-1)}{2} + 1, 3m (m-1)+m \epsilon \right),\label{eq:kv}\end{align} 
then there exists a non-constant morphism of curves $f_1 \colon X \rightarrow Z$ such that $Y$ has minimum density degree $d/\deg f_1$.
\end{theorem}
The leading term in \cref{eq:kv} is $3/4d^2$, so the bound in \Cref{cor:kadets_vogt} is roughly a factor $16/3$ worse than the bound in \Cref{thm:kadets_vogt}. However, in return for this we get a generalisation in two different directions. First of all, \Cref{thm:kadets_vogt} only gives a source of degree $d$ points. However, it does not guarantee that all but finitely many closed points of degree $d$ come from this source. That is, it could very well be that there are infinitely many closed degree $d$ points on $X$ that do not come from pulling back degree $d/ \deg f_1$ points along the function $f_1$ in \cref{thm:kadets_vogt}. Secondly, the integer $d$ in \Cref{thm:kadets_vogt} is restricted to be the minimum density degree, while in \Cref{cor:kadets_vogt} the integer $d$ is allowed to be any integer that is small enough. This allows for the study of closed points of degree larger than the minimum density degree. This second generalisation is particularly relevant in light of \cite{DHJOdensitydegree} where all modular curves $X_0(N)$ of minimum density degree $5$ were determined. There, the authors of that article were unable to determine all the $X_0(N)$ with infinitely many degree $5$ points. For the 30 values of $N$ where they were unable to determine the infinitude of the degree $5$ points on $X_0(N)$ the curve $X_0(N)$ has minimum density degree $< 5$. And this causes serious difficulty due to the lack of tools to study the degree $d$ points beyond the minimum density degree. The genus bound in \Cref{cor:kadets_vogt} is too large to be of any help in these 30 cases. However, I expect the bound from \Cref{cor:kadets_vogt} to be far from sharp. So, if this bound can be reduced enough, then it might be helpful in studying the degree $5$ points on $X_0(N)$. See \cite[\S 7.1]{DHJOdensitydegree} for more details on the difficulties with respect to the 30 remaining cases.

\subsection{Comparing \texorpdfstring{\Cref{cor:kawhaja_siksek}}{Corollary 1.4}}\label{sec:single_source_comp}

\Cref{cor:kawhaja_siksek1} of \Cref{cor:kawhaja_siksek} is a generalisation of the recent Single Source Theorem by Kawhaja and Siksek.

\begin{theorem}[{Single Source Theorem \cite[Thm. 1.1]{KhawajaSiksekSingleSource}}]\label{thm:single}
Let $X$ be a nice curve over a number field $K$,
and write $J$ for the Jacobian of $X$.
Let $c \ge 2$ and suppose 
\begin{equation}\label{eqn:glb}
	\begin{cases} \genus(X)>(c-1)^2 & \text{if $c \ge 3$}\\
		\genus(X)\ge 3 & \text{if $c=2$}.
	\end{cases}
\end{equation}
Suppose that $A(K)$ is finite for every abelian subvariety $A$ of $J$ of dimension $\leq c/2$.
If $X$ has infinitely many primitive points of degree $c$,
then there is a degree $c$ morphism $f: X \rightarrow \PP^1$
defined over $K$ such that all but finitely many primitive points of degree $c$ are of the form $f^{*}(x)$ for some $x \in Y(K)$.
\end{theorem}
Note that the requirement for the genus in \Cref{cor:kawhaja_siksek} is asymptotically a factor 4 larger. However, in exchange for that the finiteness condition for every abelian subvariety $A$ of $J$ of dimension $\leq c/2$ can be dropped. This is very useful since this means that \Cref{cor:kawhaja_siksek} is applicable to any curve with infinitely many primitive points and of large enough genus irregardless of the structure of it's Jacobian.

Also note that in \Cref{cor:kawhaja_siksek} all the primitive points always come from a function of degree $c$ and never from a map to a positive-rank elliptic curve. The reason for this is that their assumptions imply that $J$ does not contain a positive rank elliptic curve. Hence, they do not need to consider this possibility. 

Finally, the Single Source Theorem due to Kawhaja and Siksek only concerns primitive points of a single degree. However, as \Cref{cor:kawhaja_siksek1} of \Cref{cor:kawhaja_siksek} shows, something stronger is true. Namely, if a curve has infinitely many primitive points of small degree, then all but finitely many of these points have to actually have to be concentrated in a single degree. This shows that there cannot be two single sources of small-degree primitive points of different degrees. This is something that is not ruled out by the results of Kawhaja and Siksek.

\section{Background}

Let $X$ be a nice curve over a number field $K$. With a closed points $x \in X$ we mean point on $x$ that is closed in the Zarisky topolgy on $X$. The residue field $K(x)$ of such a point will be a finite extension of $K$ and the degree of such a point is the degree of the field extension $[K(x):K]$. There is a bijection between the closed points of $X$ and the Galois orbits of points in $X(\overline K)$. Under this bijection, the degree of the point corresponds to the size of the Galois orbit. If $x$ is a closed point then $[x] \in \Pic^{d} X(K)$ denotes the isomorphism class of the line bundle associated to $x$. For more details on closed points and Galois orbits see \cite[\S 2.1]{VirayVogt}.

\begin{definition}[$\mathbb{P}^1$-parameterised point]
    Let $X$ be a nice curve over a number field $K$. A \emph{$\mathbb{P}^1$-parameterised point} of degree $d$ on $X$ is a closed point $x \in X$ such that there exists a non-constant morphism $f: X \to \mathbb{P}^1_K$ of degree $d$ such that $f(x) \in \mathbb{P}^1_K(K)$
\end{definition}

We let $\Sym^d X$ denote the d-th symmetric power. The points on $\Sym^d X$ correspond to effective divisors of degree $d$ on $X$, and as such there is a map $\phi :\Sym^d X \to \Pic^d X$. The scheme0theoretic image of $\phi$ will be denoted by $W_d X$. For a closed point $x$ on has that $[x]$ is contained in $W_d X (K)$. 

\begin{definition}[AV-parameterised point]
    Let $X$ be a nice curve over a number field $K$. A closed point $x \in X$ of degree is called \emph{AV-parameterised} if there exists a positive rank abelian variety $A \subseteq \Pic^0 X$ such that $[x] +A \subseteq W_d$.
\end{definition}

With these definitions, one can summarise the main theorem of \cite{FaltingsAV} as follows:
\begin{theorem}[{\cite[p. 1]{FaltingsAV}}]\label{thm:faltings}
Let $X$ be a nice curve over a number field $K$ and $d$ be a positive integer. Then all but finitely many closed points of degree $d$ on $X$ are either $\P^1$-parameterised or AV-parameterised.
\end{theorem}

The above Theorem due to Faltings is the key ingredient of the proof of the main Theorem.

With a little extra work one can also show that the existence of a single $\P^1$-parameterised or AV-parameterised point of degree $d$ implies the existence of infinitely many closed points of degree $d$, see \cite[Thm. 4.2]{belov}. 

The proof of the main Theorem also uses an elegant application of the Castelnuovo-Severi inequality.

\begin{proposition}[Castelnuovo-Severi inequality]\label{tm:CS}
Let $K$ be a perfect field and let $X,\ Y, \ W$ be curves over $K$. Let $f:X\rightarrow Y$ and $g:X\rightarrow W$ be non-constant morphisms of degree $m$ and $n$, respectively. Assume that there is no morphism $X\rightarrow Z$ of degree $>1$ through which both $f$ and $g$ factor. Then the following inequality holds:
\begin{equation} \label{eq:CS}
\genus(X)\leq m \cdot \genus(Y)+n\cdot \genus(Z) +(m-1)(n-1).
\end{equation}
\end{proposition}

See \cite[Theorem 3.11.3]{Stichtenoth09} for a proof of this inequality.

\section{Proof of \texorpdfstring{\Cref{thm:main}}{Theorem 1.1}}

By \Cref{thm:faltings}, all but finitely many closed points of degree $d$ on $X$ are either $\P^1$-parameterised or $AV$-parameterised, so it suffices to restrict to these parameterised points. Note that \cref{eq:cs_ineq} forces $\deg f > 1$. So, the theorem follows from \Cref{prop:p1_case,prop:av_case}.

The following lemma is a variation of \cite[Prop. 1.ii]{frey1994curves}
\begin{lemma}\label{lem:frey}
Let $X$ be a curve over a number field $K$, $D$ be an effective divisor of degree $d$, and $A \subset \Pic^0(X)$ a positive-dimensional abelian variety such that $A + D \subseteq W_d$. Then $\dim_K H^0(X,\O_X(2D)) > 1$.
\end{lemma}
\begin{proof}
If $L$ is a field extension of $K$ then $$\dim_L H^0(X_L,\O_{X_L}(2D)) = \dim_K H^0(X,\O_X(2D)).$$ So after replacing $K$ with a finite extension of $K$ one can assume $X(K) \neq \emptyset$ and $\# A(K) = \infty$.
Let $P \in A(K)$ be a point of infinite order. Because $X(K) \neq \emptyset$ one can find a divisor $F$ of degree 0 such that $P=\O_X(F)$. For every $i \in \Z$ choose $D_i$ be an effective divisor of degree $d$ that is linearly equivalent to $iF+D$, this is possible since  $A + D \subseteq W_d$. Now consider the divisors of the form $D_i+D_{-i} \in \Sym^{2d} X$. Then there is an integer $i$ such that $D_i+D_{-i} \neq 2D$, the reason for this is that there are only finitely many distinct ways to write $2D$ as a sum of two effective divisors. However, the linear equivalence $D_i+D_{-i} \sim ~ iF+D+-iF+D = 2D$ implies that there exists some non-constant function $f$ such that $\ddiv f = D_i+D_{-i} - 2D$. The lemma follows since  $H^0(X,\O_X(2D))$ contains the non constant function $f$.
\end{proof}

\begin{corollary}\label{cor:frey}
In the same setting as \Cref{lem:frey}, if we additionally assume that $D=x$ for some closed point $x \in X$, then either there is a function $f: X \to \P^1$ such $f^*(\infty)=x$ or there is an $f$ such that $f^*(\infty) =2x$.
\end{corollary}

\begin{proposition}\label{prop:p1_case}
Let $X$ and $Y$ be nice curves over a number field $K$, and $f: X \to Y$ be a dominant morphism of degree $>1$. Let $d$ be an integer such that 
\begin{align}
   \genus(X) > \deg(f) \genus(Y) + (d-1)(\deg f-1) \label{eq:cs_ineq_p1}.
\end{align}
Finally, let $x \in X$ be a $\P^1$-parameterised point of degree $d$. Then there exists a decomposition $X \stackrel{f_1}{\to} Z \stackrel{f_2}{\to} Y$ of $f$ with $\deg f_1 > 1$ such that $x$ comes from pulling back a closed point of degree  $d/\deg f_1$ on $Z$.
\end{proposition}
\begin{proof}
By definition of $\P^1$-parameterised $x = h^*(\infty)$ for some $h: X \to \P^1_K$ of degree $d$. Due to \cref{eq:cs_ineq_p1}, the Castelnuovo-Severi inequality \Cref{tm:CS} for $f$ and $h$ is not satisfied. In particular, there is some map of curves $f_1: X \to Z$ such that $f$ and $h$ both factor through $f_1$
\begin{figure}[h]
    \centering

\begin{tikzcd}
& X\arrow[ldd,"g",swap] \arrow[d,"f_1"] \arrow[rdd,"f"] & \\
& Z  \arrow[ld,"h"] \arrow[rd,"f_2",swap]& \\
\P^1 & & Y 
\end{tikzcd}
    \caption{}
    \label{fig:cs_diagram}
\end{figure} Let $h: Z \to \P^1$ be such that $g = h \circ f_1$, then $x=g^*(\infty) = f_1^*(h^*(\infty))$. In particular $x$ is the pullback along $f_1$ of $h^*(\infty)$ which is of degree $\deg h = d/\deg f_1$.

\end{proof}

\begin{proposition}\label{prop:av_case}
Let $X$ and $Y$ be nice curves over a number field $K$, and $f: X \to Y$ be a dominant morphism of degree $>1$. Let $d$ be an integer such that 
\begin{align}
   \genus(X) > \deg(f) \genus(Y) + (2d-1)(\deg f-1) \label{eq:cs_ineq_av}.
\end{align}
Finally, let $x \in X$ be an $AV$-parameterised point of degree $d$ at which $f$ is unramified. Then there exists a decomposition $X \stackrel{f_1}{\to} Z \stackrel{f_2}{\to} Y$ of $f$ with $\deg f_1 > 1$ such that $x$ comes from pulling back a closed point of degree $d/\deg f_1$ on $Z$.
\end{proposition}

\begin{proof}
Since $x$ is $AV$-parameterised the conditions of \Cref{lem:frey} are satisfied so that \Cref{cor:frey} can be applied to get an $g: X \to \P^1$ such that either $g^*(\infty)=x$ or $g^*(\infty)=2x$. In the first case $x$ is $\P^1$ parameterised as wel so that \Cref{prop:p1_case} can be applied.

Now consider the case $g^*(\infty)=2x$, so that $\deg g=2d$. By \cref{eq:cs_ineq_av} and the Castelnuovo-Severi inequality, again there has to be a diagram as in \cref{fig:cs_diagram} and hence there is a decomposition $X \stackrel{f_1}{\to} Z \stackrel{h}{\to} \P^1$ of $g$. This time $2x=g^*(x)=f_1^*(h^*(\infty))$. Since $f_1$ is unramified at $x$ and $2x=f_1^*(h^*(\infty))$ it follows that $h^*(\infty)=2y$ for some closed point $y$ on $Z$. For this closed point $y$ one has $f_1^*(y)=x$ and $\deg y = \frac 1 2 \deg h = \frac 1 2 \deg g /\deg f_1=d/\deg f_1$, and the theorem follows.
\end{proof}

\section{Proof of \texorpdfstring{\Cref{cor:kawhaja_siksek}}{Corollary 1.4}}
\Cref{cor:kawhaja_siksek2}:
The proof starts by describing how to obtain the map $f: X\to Y$. Let $f' : X \to \P^1_K$ be a non-constant morphism of minimal degree. Since $X$ has infinitely many closed degree $c$ points, the map $f'$ has degree at most $2c$ \cite[Prop. 2]{frey1994curves}. For every decomposition $X \stackrel{f'_1}{\to} Z \stackrel{f'_2}{\to} \P^1_K$ with $\deg f_1>1$ of $f'$ we chose one representative up to isomorphism and let $S'$ be the set of all $f_1$ coming from these representatives. Applying \Cref{thm:main} to $f'$, one concludes that all but finitely many closed degree $c$ points come from pulling back a degree $c/ \deg f_1$ point along some $f_1 \in S$. In particular, since $S$ is finite, one can chose an $(f: X\to Y) \in S$ such that infinitely many primitive points of degree $c$ come from pulling back a closed degree $c/\deg f$ point on $Y$.

Now for this choice of $f: X\to Y$ one has $\deg f=c$ and $\#Y(K) =\infty$. To see this, let $x$ be a primitive point of degree $c$, which is the pullback of a closed point $z$ on $Z$ of degree $c/\deg f$. Since $x$ is primitive and we have inclusions $K \subseteq K(z) \subseteq K(x)$. It follows that $K(z) = K$; $K(z) = K(x)$ cannot happen since $[K(x):K(z)]= \deg f>1$. In particular, $\deg f=c$ and the infinitely many primitive degree $c$ points on that are pullbacks along $f$, are pullbacks of $K$ rational points on $Y$ so that $\#Y(K) = \infty$. In particular, $Y$ has genus $0$ or $1$. 

The morphism $f$ is indecomposable, since a decomposition of $f$ would lead to a field intermediate to $K=K(z)$ and $K(x)$, which does not exist by the primitivity of $x$.

Observe that $$\genus(X) > (2c-1)^2 \geq c+(2c-1)(c-1) \geq \deg f\genus(Y)+(2c-1)(\deg f -1).$$ So the fact that all but finitely many primitive points of degree $c$ are of the form $f^{*}(x)$ for some $x \in Y(K)$ follows from applying \Cref{thm:main} to $f$, and using the fact that $f$ is indecomposable. 

\Cref{cor:kawhaja_siksek4}: One easily deduces the inequality $$\genus(X)>2d(c-1)+1=c+(2d-1)(c-1)\geq \genus(Y)\deg f+(2d-1)(\deg f-1).$$ Therefore, this follows from applying \Cref{thm:main} to $f$, and using the fact that $f$ is indecomposable.

\Cref{cor:kawhaja_siksek3}: This follows immediately from \Cref{cor:kawhaja_siksek4}.

\Cref{cor:kawhaja_siksek1}:  This follows immediately from \Cref{cor:kawhaja_siksek4}. Indeed, all but finitely many points on $X$ of degree $d$ are of the form $f^*(z)$ where $z$ is a closed point of degree$d/c$ on $Z$. For such a point to be primitive, one needs that $z$ has degree $1$. This is impossible since $c\neq d$.

\printbibliography
\end{document}